\documentclass[11pt]{article}

\usepackage{amsfonts}
\usepackage{amssymb}
\usepackage{graphicx}

\usepackage{ifthen}
\usepackage{color}
\usepackage{cite}
\input{steve.lib}

\begin{document}

\title{Request-Based Gossiping without Deadlocks\thanks{
%J. Liu and A. S. Morse are with Yale University, USA
%(\texttt{\{ji.liu, as.morse\}@yale.edu}).
Preliminary versions of this paper have appeared in the Proceedings of the
50th and 51st IEEE Conference on Decision and Control \cite{cdc11req,cdc12gossip}.
J. Liu is with the Coordinated Science Laboratory, University of Illinois at Urbana-Champaign, USA (\texttt{jiliu@illinois.edu}).
S. Mou is with Purdue University, USA (\texttt{mous@purdue.edu}).
A. S. Morse is with Yale University, USA (\texttt{as.morse@yale.edu}).
B. D. O. Anderson and C. Yu are with the Australian National University
and National ICT Australia Ltd., Australia
(\texttt{\{brian.anderson, brad.yu\}@anu.edu.au}).
C. Yu is also with Shandong Computer Science Center, Jinan, China.
}}

\author{J. Liu  \hspace{.2in}  S. Mou  \hspace{.2in}  A. S. Morse  \hspace{.2in}  B. D. O. Anderson  \hspace{.2in}  C. Yu}
\maketitle

\begin{abstract}
By the distributed averaging problem is meant the problem of computing the average value of a
set of numbers possessed by the agents  in a distributed network using only  communication between neighboring agents.
Gossiping is a well-known approach to the  problem which seeks to iteratively arrive at a solution  by
allowing each agent to interchange information with at most one neighbor  at each iterative step. Crafting a
gossiping protocol which accomplishes
this is challenging because gossiping is an inherently collaborative process which
can lead to deadlocks unless careful precautions are taken to ensure that it does not.
Many gossiping protocols are request-based which means simply
that a gossip between two agents will occur whenever one of
the two agents accepts a request to gossip placed by the other.
In this paper, we present three deterministic request-based protocols.
We show by example that the first can deadlock.
The second
is guaranteed to avoid deadlocks and  requires fewer transmissions per iteration than
standard broadcast-based distributed averaging protocols
by exploiting the idea of local ordering together with the notion of an agent's neighbor queue;
the protocol requires
the simplest queue updates, which provides an in-depth understanding of how local
ordering and queue updates avoid deadlocks.
 It is shown that a third protocol
which uses a slightly more complicated queue update rule can lead to significantly
faster convergence; a worst case bound on convergence rate is provided.
%An asynchronous version of the third protocol is also provided
%which is shown to solve the distributed averaging problem in the case when
%each agent has its own clock.
\end{abstract}

%\noindent{\bf Key words:} gossiping, cooperative control, multi-agent systems, convergence rates, asynchronous systems

\section{Introduction}

Over the past decade, there has been considerable interest in
developing algorithms for distributed computation and decision making
among the members of a group of sensors or mobile autonomous agents via local interactions.
Probably the most notable among these  are those algorithms intended to
 cause such a  group  to reach
a consensus  in a distributed manner
\cite{Ts3,vicsekmodel,reza1,luc,ReBe05,HeTs13,cdc14}.
Consensus processes have found applications in a wide range of fields including
%robotic teams \cite{beard},
web search engines \cite{tempo},
social networks \cite{cdc12}, and electric power grids \cite{bullo}.
Many variants of consensus problems, such as %asynchronous consensus \cite{ming3},
quantized consensus \cite{basar}, constrained consensus \cite{nedic},
and modulus consensus \cite{modulus}, have been proposed lately.

We are interested in distributed averaging,
a particular type of consensus process which has received much attention recently \cite{boyd04}.
A typical distributed averaging process deals with a network of $n>1$ {\em agents} and
 the constraint that
  each agent $i$ is able to communicate
 only with  certain other agents called agent $i$'s {\em neighbors.} Neighbor relationships are
conveniently characterized by a simple, undirected, connected graph $\bbb{A}$  in which vertices correspond to agents
 and edges indicate neighbor relationships. Thus   the
  neighbors  of an agent $i$ have the same labels as the  vertices in $\bbb{A}$
   which are adjacent to vertex $i$.
  Initially, each agent $i$
has  or  acquires    a real number $y_i$ which might be  a measured temperature or something similar.
 The  {\em distributed averaging problem}  is to devise
 an algorithm which will enable each agent
to compute the average
$$y_{{\rm avg}} = \frac{1}{n}\sum_{i=1}^ny_i$$
using information received only from  its neighbors.
%The problem is also sometimes called the average consensus problem \cite{chamie}.
There are many variants of the problem.
 For example, the edges of $\bbb{A}$ may change over time \cite{metro2}.
 Another variant  assumes that communications between neighboring agents can be unidirectional \cite{acc12}.
This paper  considers the case when $\bbb{A}$ does not depend on time
and communications between neighboring agents are all bidirectional.

There are three important approaches to the distributed averaging problem
in the literature: linear iterations \cite{boyd04},
gossiping \cite{boyd052}, and double linear iterations \cite{acc12}
(which are also known as
push-sum algorithms \cite{kempe}, weighted gossip \cite{weighted}, and ratio consensus \cite{wcdc12}).
Double linear iterations are specifically tailored to the case in which unidirectional communications exist;
they can solve the problem when $\bbb{A}$ is directed, strongly connected,
but under the assumption that each agent is aware of the number of its
out-going neighbors.
Both linear iterations and gossiping work for the case in which
all communications between neighbors are bidirectional; in this case,
double linear iterations have the disadvantage that they require
updating and transmission of an additional variable for each agent.

Linear iterations are a well studied approach to the problem
in which  each agent communicates with all of its neighbors on each iteration,
and thus are sometimes called {\em broadcast} algorithms.
It is clear that broadcast algorithms typically require a lot of transmissions between neighbors per unit time,
which may not be possible to secure in some applications,
particularly when communication cost is an important issue on each iteration.
For example, fewer transmissions per iteration can increase the time interval between any two
successive recharges of a sensor, and improve the security of the
network by reducing the opportunities of being hacked or eavesdropped.

Gossiping is an alternative approach to the distributed averaging problem which does not involve broadcasting.
%The idea of gossiping is very simple.
%A pair of neighbors with labels $i$ and $j$ are said to {\em gossip} at time $t$ if  both $x_i(t+1)$ and $x_j(t+1)$
%are set equal to the average of $x_i(t)$ and $x_j(t)$.
An important rule of gossiping is that
each agent is allowed to gossip with at most one neighbor at one time.
This is the reason why gossiping algorithms  do not involve
broadcasting. Thus gossiping algorithms have the potential
to require less transmissions  per iteration than broadcast algorithms.
Moreover, the peer-to-peer nature of gossiping simplifies the implementation of
algorithms and reduces computation complexity on each agent.
As a trade-off, one would not expect gossiping algorithms to
converge as fast as broadcast algorithms.

Most existing gossiping algorithms are probabilistic in the sense that
the actual sequence
  of gossip pairs which occurs
  during a specific gossip process is determined  probabilistically \cite{boyd052,fabio}.
Recently, deterministic gossiping has received some attention \cite{low&murray,pieee}.
Probabilistic gossiping algorithms aim at achieving consensus asymptotically with probability one,
whereas deterministic gossiping algorithms are intended to guarantee that under all conditions,
  a consensus will  be achieved asymptotically.
Both approaches have merit.  The probabilistic approach is  easier
both in terms of algorithm
development and convergence analysis.  The deterministic approach forces one
 to  consider worst case
scenarios   and has the potential of yielding algorithms which may outperform
 those obtained using the  probabilistic approach.
For example, the deterministic approach rules out the possibility of deadlocks
which may occur in probabilistic gossiping algorithms.

   Crafting a deterministic protocol is challenging because gossiping is an inherently collaborative process which
can lead to deadlocks unless careful precautions are taken to ensure that it does not.
 The global ordering \cite{low&murray},
centralized scheduling \cite{pieee}, and broadcasting \cite{alex} are the existing ways
to avoid deadlocks.
Both global ordering and centralized
scheduling require a degree of network-wide coordination
and broadcasting requires each agent to obtain the values
of {\em all} of its neighbors' ``gossip variables'' at each clock time,
which may not be possible to secure in some applications.
The aim of this paper is to present %several
deterministic gossiping protocols which do not utilize global ordering, centralized
scheduling, or broadcasting and
are guaranteed to solve
the distributed averaging problem.

Three gossiping protocols are considered in the paper.
We show by example that the first can deadlock. After minor
modifications, a second protocol is obtained.
The second protocol
is guaranteed to avoid deadlocks and  requires fewer transmissions per iteration than
standard broadcast-based distributed averaging protocols,
%The second protocol, first presented in \cite{cdc12gossip},
which requires the simplest queue updates and thus provides an in-depth understanding of how local
ordering and queue updates avoid deadlocks.
 It is shown both by analysis and computer studies that a third protocol
which uses a slightly more complicated queue update rule can lead to significantly
faster convergence.
%All three protocols exploits the ideas of local ordering proposed in \cite{alex}
%together with the notion of an agent's ``neighbor queue''.
%An asynchronous version of the third protocol is also proposed
%which is shown to solve the distributed averaging problem even when
%each agent has its own independent clock.

The material in this paper was partially presented in \cite{cdc11req,cdc12gossip},
but this paper presents a more comprehensive treatment of the work.
Specifically, the paper %provides an asynchronous protocol in Section \ref{asyn},
provides proofs for Theorems \ref{complete2}, \ref{tree}, \ref{anu},
Propositions \ref{single}, \ref{morse},
Lemmas \ref{gossip}, \ref{pizza}, \ref{lib}, \ref{iff}, \ref{contain},
and establishes an additional result Lemma \ref{dstep}, which were
not included in \cite{cdc11req,cdc12gossip}.
Note that Protocol III in the paper was briefly outlined in \cite{pieee},
but without a proof of correctness.

\section{Gossiping}

As  mentioned in the introduction, gossiping is a way of distributed averaging.
Consider a group of $n>1$ agents
labeled $1$ to $n$.\footnote{
The purpose of labeling of the agents is only for convenience. We do
not require a global labeling of the agents in the network. We only assume that each agent can differentiate its neighbors.}
 Each agent $i$ has control over  a real-valued scalar
quantity $x_i$ called agent $i$'s {\em gossip variable} whose value $x_i(t)$ at time $t$ represents
agent $i$'s estimate of the desired average $y_{{\rm avg}}$ at that time.
A {\em gossip} between agents $i$ and $j$, written $(i,j)$,
occurs at time $t$
 if the values of both agents' variables at time $t+1$  equal
  the average of their
 values at time $t$. In other words, $x_{i}(t+1)=x_j(t+1) = \frac{1}{2}(x_i(t)+x_j(t))$.
 If agent $i$ does not gossip at time $t$, its gossip variable does not change; thus
 in this case $x_i(t+1) = x_i(t)$.
Generally not every pair of agents
  is allowed to gossip. The edges of a simple, undirected, connected graph $\bbb{A}$
specify which pairs of agents are allowed to gossip.
  In other words,
   a gossip  between agents $i$ and $j$
  is  {\em allowable} if $(i,j)$ is an edge in $\bbb{A}$.
   We  sometimes call $\bbb{A}$ an {\em allowable gossip graph}.

An important rule of gossiping is that in a gossiping process,
  each agent is allowed to gossip with at most one of its neighbors at one time.
This  rule
does not preclude the possibility
  of two  or more  pairs of agents gossiping at the same time,
  provided that the pairs have no agent in common. To be more precise, two gossip
   pairs $(i,j)$ and $(k,m)$
  are {\em noninteracting} if neither $i$ nor $j$ equals either $k$ or $m$.
When multiple noninteracting  pairs of allowable gossips  occur simultaneously,
 the simultaneous occurrence of all such gossips is called a
  {\em  multi-gossip}. In other words, a multi-gossip at time $t$ is the set
  of all gossips which occur at time $t$ such that each such pair is allowable
   and that any two of
  such pairs are noninteracting.

 Gossiping processes can be modeled
     by    a discrete-time linear system of the form
    \eq{x(t+1) = M(t)x(t),\;\;\;\;\; t=0,1,2,\ldots \label{xxx}}
    where $x\in\R^n$ is
     a state vector of gossiping variables
    and $M(t)$ is a matrix characterizing how $x$ changes as the result of the
     gossips which take place at time $t$.
%sometimes $M(t)$    depends      on $x$ although the notational dependence is often suppressed.
     %For the present it will be convenient to
     %ignore the possibility of multi-gossips in which
     %case each $M(t)$ must have a very particular structure.
 If a single pair of  agents
    $i$ and $j$ gossip
     at time $t\geq 0$, then $M(t) = P_{ij}$
      where
      $P_{ij}$ is the $n\times n$   matrix for which
       $p_{ii} = p_{ij}= p_{ji} = p_{jj} = \frac{1}{2}$,
       $p_{kk} = 1$, $k\not\in\{i,j\}$,
       and all remaining entries equal $0$. We call such $P_{ij}$ a
       {\em single gossip primitive gossip matrix}.
        For convenience, we include in the set of primitive
        gossip matrices, the $n\times n$ identity matrix $I$;
        the identity matrix can be thought of as the update matrix to
        model the case in which no gossips occur at time $t$.
 If  a multi-gossip   occurs at time $t$,
         then as a consequence of non-interaction, $M(t)$ is simply the product
         of the single gossip primitive  gossip matrices corresponding to the individual gossips
         comprising the multi-gossip; moreover, the primitive gossip matrices in the product commute with each
           other and thus
          any given permutation of the single gossip primitive matrices in the product determines
          the same matrix $P $.  We call $P$ the {\em primitive gossip matrix} determined
           by the multi-gossip  under consideration.  For example, if
             agents $a,b,$ and $c$ gossip with
          $d,e,$ and $f$ respectively at time $t$, then the   primitive gossip matrix determined
          is $P=P_{ad}P_{be}P_{cf}$ and thus in this case
           $M(t)=P$ and $x(t+1) = Px(t)$.

We will see that for any gossiping process determined by the protocols presented in this paper,
the update matrix $M(t)$ in \rep{xxx} also depends on the state $x(t)$ and thus
$$x(t+1) = M(x(t),t)x(t),\;\;\;\;\;\; t=0,1,2,\ldots $$
while each $M(x(t),t)$ is still a primitive gossip matrix. Therefore, the system to be studied is essentially nonlinear,
which is a significant difference from those in \cite{boyd052,low&murray}.
This difference also makes the protocol design and analysis more challenging than probabilistic protocols.

\subsection{Generalized Gossiping}

Although in this paper
   we shall be interested in gossiping protocols which stipulate that each agent is allowed
   to gossip with at most one of its neighbors at one time,
    as we shall see later, there is value in taking the time here to generalize the idea.

    Let us
    agree call a  subset  $\scr{L}$  of $m>1$  agent   labels, a {\em neighborhood}
    if each  pair of distinct labels in  $\scr{L}$
 are the labels of   vertices in $\bbb{A}$ which are connected.
We say that the  agents with labels in $\scr{L}$
perform a {\em gossip of order $m$}  at time $t$ if each updates its
gossip variable to the average of all; that is,
 if $x_{i}(t+1) = \frac{1}{m}\sum_{j\in\scr{L}}x_j(t),\;i\in\scr{L}$. A {\em generalized gossip}
 is a gossip of any order.
A gossip    without the modifier ``generalized'', will  continue to mean
  a gossip of order $2$.
 A {\em generalized multi-gossip} at time $t$ is a finite set of
  generalized gossips with disjoint neighborhoods which  occur simultaneously   at time  $t$.

It is worth emphasizing that the concepts of generalized gossips
and multi-gossips are introduced only for the purpose of analysis.
Generalized gossips and multi-gossips do not occur in any gossiping sequence
generated by the protocols presented in this paper.
But the effect of ``virtual gossips'' generated by the protocols in this paper
is the same as the occurrence of generalized (multi-)gossips;
see \S \ref{gseq} for detailed explanation.

The idea of a primitive gossip matrix extends naturally  to generalized  gossips.
In particular, %if $\scr{L}$ is a neighborhood % whose labels are the labels of to a connected subgraph of $\bbb{A}$,
 we associate with   a neighborhood  $\scr{L}$
the $n\times n$  doubly stochastic matrix $P_{\scr{L}}$ where
$ p_{jk} =\frac{1}{m+1}, \;j,k \in\scr{L}$,
$p_{jj} = 1,\;j\not\in \scr{L}$,  and $0$s elsewhere. We call $P_{\scr{L}}$ the {\em
primitive gossip matrix} determined by $\scr{L}$. By the {\em graph induced by} $P_{\scr{L}}$, written
 $\bbb{G}_{\scr{L}}$,
 we mean the
 spanning subgraph
 of $\bbb{A}$ whose edge set is all edges in
  $\bbb{A}$ which are incident on   vertices with labels which are both in $\scr{L}$.
More generally, if $\scr{L}_1,\scr{L}_2,\ldots,\scr{L}_k$
are $k$ disjoint neighborhoods, the matrix $P_{\scr{L}_1}P_{\scr{L}_2}\cdots P_{\scr{L}_k}$  is the
primitive gossip matrix determined by  $\scr{L}_1,\scr{L}_2,\ldots,\scr{L}_k$ and the graph induced by
 $P_{\scr{L}_1}P_{\scr{L}_2}\cdots P_{\scr{L}_k}$ is
the union
of the induced graphs $\bbb{G}_{\scr{L}_i},\;i\in\{1,2,\ldots,k\}$. Note that the matrices in the product
$P_{\scr{L}_1}P_{\scr{L}_2}\cdots P_{\scr{L}_k}$ commute because the $\scr{L}_i$ are disjoint so the order of the matrices in the product is not
important for the definition to make sense. Note also that there are only finitely many   primitive
 gossip matrices associated with $\bbb{A}$.

   \subsection{Gossiping Sequences}\label{gseq}

Let  $\gamma_1,\gamma_2,\ldots $ be an infinite sequence of multi-gossips
corresponding to some or all of the edges in  $\bbb{A}$.
Corresponding to such a  sequence
  is a sequence of  primitive gossip matrices $Q_1, Q_2,\ldots $ where $Q_i$ is the primitive
   gossip matrices of the  $i$th multi-gossip in the sequence.
For given $x(0)$, such a gossiping  matrix sequence  {\em generates}
 the sequence of vectors
\eq{ x(t) = Q_{t}Q_{t-1}\cdots Q_{1}x(0),\;\;\;\;\;\;t>0\label{pp}}
which we call a {\em gossiping  sequence}.    %We call any  protocol \{or algorithm\}
 %which generates such a sequence, a {\em gossip protocol} \{or algorithm\}.
  We have purposely restricted this definition of a gossiping sequence to multi-gossip sequences, as opposed to
{\em generalized} multi-gossip sequences, since we will only be dealing with algorithms involving multi-gossips.
Our reason for considering  generalized multi-gossips  will become clear in a moment.

As will soon be obvious,  the matrices $Q_i$ in \rep{pp}
 are not necessarily the only  primitive gossip matrices
  for which \rep{pp} holds.  This non-uniqueness
 can play a crucial rule in understanding certain gossip protocols which are not  linear iterations.
 To understand why this is so,  let us agree to say that the transition $x(\tau)\longmapsto x(\tau+1)$ contains a
 %$mth$ order
{\em virtual gossip} if there is  a neighborhood   $\scr{L}$  % with $m$ labels
for which
$x_i(\tau)=x_j(\tau),\;i,j\in\scr{L}$.  We say that agent $i$ has gossiped virtually with agent $j$ at time $t$,
if $i$ and $j$
 are both labels in $\scr{L}$.  Thus while   we are only interested in algorithms in which an agent may
  gossip with at most one neighbor at any one time,  for such algorithms there may be times at which
   virtual gossips  occur between an agent and  one or more  of its neighbors.
Suppose that  for some time  $\tau< t$,  the transition $x(\tau)\longmapsto x(\tau+1)$
contains such a
virtual gossip   and let $P_{\scr{L}}$  denote the  primitive gossip matrix determined
by $\scr{L}$.
Then clearly $P_{\scr{L}}x(\tau) = x(\tau)$ which means that the matrix
$Q_{\tau+1}$  in the product $Q_{t}Q_{t-1}\cdots Q_{1}$ can be replaced by the matrix $Q_{\tau+1}P_{\scr{L}}$
without changing the validity
of \rep{pp}. Moreover $Q_{\tau+1}P_{\scr{L}}$ will be a primitive gossip matrix if the neighborhoods which define
$Q_{\tau+1}$  are disjoint with $\scr{L}$.
 The importance of this elementary observation is simply this. Without taking into account
virtual gossips in equations such as \rep{pp},  it may  in some cases to be impossible to conclude
that the matrix  product  $Q_{t}Q_{t-1}\cdots Q_{1}$ converges as $t\rightarrow \infty$ even though
the gossip  sequence  $x(1),x(2), \ldots $ does. Later in this paper we will describe a gossip
protocol for which this is true.

Prompted by the preceding, let us agree to say that
  a gossiping sequence satisfying \rep{pp} is {\em consistent} with a
sequence of  primitive gossip matrices
$P_1,P_2,\ldots $ if
\eq{ x(t) = P_{t}P_{t-1}\cdots P_{1}x(0),\;\;\;\;\;\;t>0\label{pp2}}
It is obvious that if the sequence $x(t),\;t\geq 0$ is consistent with the sequence
 $P_1,P_2,\ldots $ and the latter converges,
 then so does the former. Given a gossip vector sequence,  our task then is
   to find, if possible, a consistent,  primitive gossip  matrix sequence which is also convergent.

%We will be concerned exclusively with generalized gossips which are either gossips between agent pairs
%or virtual gossips.
%By a {\em generalized gossip sequence} we mean any sequence of generalized
 %multi-gossips.
 As we have already noted,
 $\bbb{A}$  has associated with it a finite family of   primitive  gossip matrices
  and each primitive gossip matrix
 induces a spanning subgraph of $\bbb{A}$. It follows that any  finite sequence of  primitive gossip
  matrix $P_1,P_2, \ldots, P_k$  {\em induces}
 a spanning subgraph of $\bbb{A}$  whose edge set is the union of the edge sets of the graphs induced by
 all of the $P_i$.
    We say that the primitive gossip matrix  sequence
$P_1,P_2, \ldots, P_k$
  is {\em complete}
 if the graph the sequence  induces
  is a {\em connected} spanning subgraph
  of
  $\bbb{A}$.
An infinite sequence of  primitive gossip matrices
$P_1,P_2,\ldots $ is  {\em repetitively complete}    with period $T$, if each successive subsequence of
  length $T$
   in the sequence is complete. A gossiping sequence $x(t),\;t>0$
 is  {\em repetitively complete}    with period $T$, if  there is a consistent
sequence of  primitive gossip matrices which is repetitively complete with period  $T$.
The importance of  repetitive completeness is as follows.

\begin{theorem} Suppose $P_1,P_2,\dots $ is an infinite sequence of primitive gossip
matrices which is repetitively complete with period $T$. There exists a real
nonnegative number $\lambda <1$,   depending only on $T$ and  the $P_t$, for which
$$\lim_{t\rightarrow\infty}P_tP_{t-1}\cdots P_1x(0)  = y_{{\rm avg}}\mathbf{1}$$
as fast as $\lambda^t$ converges to zero.
\label{tm11}\end{theorem}
A proof of this theorem can be found in \cite{pieee}.
There are also several different ways to prove this theorem using ideas from
  \cite{vicsekmodel,luc,ReBe05}.

\section{Request-Based Gossiping}

Request-based gossiping is a gossiping process in which a gossip occurs between two agents
whenever one of the two   accepts a request to gossip placed by  the other.
The aim of this section is to design deterministic request-based gossiping protocols
which can solve the distributed averaging problem.
The design of such deterministic protocols is more complicated than probabilistic ones
since a deterministic protocol must rule out the possibility of deadlocks whereas
in a probabilistic protocol, deadlocks are allowed to occur as long as their probability
goes to zero as time goes to infinity.

%In a request-based gossiping process, a given agent
% $i$ may gossip with one of its  neighbors at time $t$ only if $t$
%is either an ``event time'' of agent $i$
%or an ``event time'' of the neighbor which has made a request to gossip with agent $i$.
%By an {\em event time} of agent $i$ is meant a time at which agent $i$ may place a
%request to gossip with  one of its neighbors. By an {\em event time interval} of agent $i$ is
% meant the interval
%of time between two successive event times of agent $i$.
%For obvious reasons, we assume that  the lengths of
% agent $i$'s event time intervals are all
% bounded above by a finite positive number $T_i$.
% We write $\scr{T}_i$ for the  set of event times of agent $i$ and
%  $\scr{T}$ for the union of the event time sequences of all $n$
%agents.

In the cases when an  agent
 who has placed a  request to gossip, at the same time receives a request to gossip from another agent,
 conflicts leading to deadlocks  can arise.
It is challenging to devise deterministic protocols which resolve such conflicts while
at the same time ensuring exponential convergence
of the gossiping process generated by the protocols.
%One way to avoid such conflicts is to assign event times  off line so that no agent can
%receive a request to gossip at any of its own event times. There are several ways
%to do this which have been discussed in \cite{pieee}.

From time to time, an agent may have more than one neighbor
to which it is able to make a request to
gossip with. Also from time to time, an agent may
receive more than one request to gossip from its neighbors.  While
in such situations decisions
about who to place a request with or whose request to accept can be randomized,
in this paper we will examine
only completely deterministic strategies.  To do this we will assume that
 each agent orders all its neighbors
 according to some priorities so when a choice occurs
among neighbors, the agent will always
choose the one with highest priority.
%Consider first the situation when the event times of each agent and each agent's
%neighbor priorities are chosen off line
%and are fixed throughout the gossiping process. Assume that the event times are chosen so that no agent
%can receive a request to gossip at any of its own event times.
% Our aim is to show that this arrangement can be problematic.
%The following  protocol illustrates this.
The simple example in \cite{pieee} illustrates that fixed priorities can
be problematic (see Protocol I in \cite{pieee} and the example which follows).
The global ordering \cite{low&murray} and centralized scheduling \cite{pieee} are the
two ways in the literature to overcome them.
Both global ordering and centralized
scheduling require certain degree of network-wide coordination which may not be possible to
secure in some applications.
  In what follows  we take an alternative approach which is fully distributed.

In the light of Theorem \ref{tm11},  we are interested in devising  gossiping protocols
 which generate repetitively complete gossip sequences. Towards this end, let us agree to say that
an  agent $i$ has completed a {\em round}
 of gossiping after it has gossiped with each of its neighbor at least once.  Thus the finite
 sequence of primitive gossiping matrices corresponding to a finite sequence of multi-gossips
 for the entire group of $n$ agents which has occurred over an interval of length $T$,
  will be  complete
 if each agent in the group completes a round of gossiping over the same interval.
% In fact
% in the case when $\bbb{A}$ is a tree,  the only way such a sequence could be
%  complete is if  over the same period
%  each agent in the group completes a round.
In \S \ref{correct}, the concept of a round of gossiping will be generalized by
taking into account virtual gossips.

For the protocols which follow it will be necessary for each agent $i$ to keep track
of where it  is in a particular round.  To do this,  agent $i$
makes use of a recursively updated     {\em  neighbor queue} $\mathbf{q}_i(t)$
 where $\mathbf{q}_i(\cdot )$ is a function from
$\scr{T}$ to the set of all possible  lists of the $n_i$ labels in
 $\scr{N}_i$, the neighbor set of agent $i$.
  Roughly speaking,  $\mathbf{q}_i(t)$ is a list of the labels of
the neighbors of
  agent $i$  which defines the queue of neighbors at time $t$ which are in line to gossip
  with agent $i$. %The updating of $\mathbf{q}_i(t)$ is straightforward: If neighbor $j$ gossips with agent $i$ at time
%  $t$, the updated queue $\mathbf{q}_i(t+1)$ is obtained by moving agent $j$'s label from its current
%   position in $\mathbf{q}_i(t)$, to the end of the queue.  If on the other hand, agent $i$ does not gossip
%    at time $t$, $\mathbf{q}_i(t+1)=\mathbf{q}_i(t)$.

%A disadvantage of Protocol III is that it requires the distinct neighbor event times assumption.
% This assumption can only be satisfied by off-line assignment of event times  for each agent, and in some applications
% such an off-line assignment may be undesirable.
   In a recent doctoral thesis \cite{alex}, a clever gossiping protocol is proposed which
 does not require the distinct neighbor event times assumption. The protocol avoids deadlocks and
 achieves consensus   exponentially  fast.
A disadvantage of  the protocol  in \cite{alex} is that it requires each agent to obtain the values of {\em all}
   of its neighbors' gossip variables at each clock time.
By exploiting one of  the key ideas  in \cite{alex} together with the
   notion of an agent's neighbor  queue $\mathbf{q}_i(t)$ defined earlier, it is possible
   to  obtain  a  gossiping protocol which also  avoids deadlocks and achieves consensus exponentially
    fast but  without requiring each agent to obtain  the values of  all of its  neighbors' gossip
     variables at each iteration.

In the sequel, we will outline a gossiping algorithm in which at time $t$, each agent $i$
 has a  single {\em preferred neighbor} whose label  $i^*(t)$ is in the front of queue $\mathbf{q}_i(t)$.
At time $t$ each agent $i$ transmits to its preferred neighbor its label $i$ and the
 current value of its gossip variable $x_i(t)$. Agent $i$ then transmits the current value of its
 gossip variable to those agents  which have agent $i$ as their preferred neighbor;  these neighbors plus
 neighbor $i^*(t)$
 are agent $i$'s   {\em receivers}  at time $t$. They are the neighbors of agent $i$ who know the
 current gossip  value of agent $i$.
 Agent $i$
 is presumed  to have placed a request to gossip with its preferred neighbor  $i^*(t)$
 if   $x_i(t)>x_{i^*(t)}(t)$;
  agent $i$ is    a
{\em requester} of
  agent $i^*(t)$ whenever this is so.   Note that while an agent has exactly one preferred neighbor,
  it may at the same time
   have anywhere from zero to $n_i$ requesters, where $n_i$ is the number of neighbors of agent $i$.

\subsection{A Raw Model}

\noindent {\bf Protocol I: }
Between clock times $t$ and $t+1$ each agent $i$ performs the steps enumerated below in the order
indicated. Although the agents' actions need not be precisely synchronized, it is understood that
for each $k\in\{1,2,3\}$ all agents complete step $k$ before any embark on step $k+1$.

\begin{enumerate}
\item {\bf 1st Transmission:} Agent $i$ sends its gossip variable value $x_i(t)$ to its current
preferred neighbor. At the same time agent $i$ receives the gossip
values from all of those neighbors which have agent $i$ as their current preferred neighbor.

\item {\bf 2nd Transmission:} Agent $i$ sends its current gossip value $x_i(t)$ to those neighbors which
have agent $i$ as their current preferred neighbor.

\item {\bf Acceptances:}

\begin{description}
  \item[(a)] If agent $i$ has not placed a request to gossip but has received at least one request
to gossip, then agent $i$ sends an acceptance to that particular requesting neighbor whose
label is closest to the front of the queue $\mathbf{q}_i(t)$.

  \item[(b)] If agent $i$ has either placed a request to gossip or has not received any requests to gossip,
then agent $i$ does not send out an acceptance.
\end{description}

\item {\bf Gossip variable and queue updates:}
\begin{description}
\item[(a)] If agent $i$ sends an acceptance to or receives an acceptance from neighbor $j$,
then agent $i$ gossips with neighbor $j$ by setting
  $$x_i(t+1) = \frac{x_i(t)+x_j(t)}{2}$$
  Agent $i$ updates its queue by moving $j$ from its current
positions in $\mathbf{q}_i(t)$ to the end of the queue.

  \item[(b)] If agent $i$ has not sent out an acceptance nor received one, then agent $i$ does not update
the value of $x_i(t)$. In addition, $\mathbf{q}_i(t)$ is not updated except when agent $i$'s gossip
value equals that of its current preferred neighbor. In this special case agent $i$
moves the label $i^*(t)$ from the front to the end of the queue.
\end{description}

\end{enumerate}

It is possible show that this protocol ensures that at each time $t$, either
$x_i(t)=x_{i^*(t)}(t)$ for some agent $i$ or a gossip must take place between two agents
whose gossip variables have different values.
But the example in \cite{cdc12gossip}
shows that this strategy will not necessarily lead to a consensus
(see Section III in \cite{cdc12gossip}).

\subsection{A Corrected Protocol}\label{correct}

It is possible to guarantee an exponentially fast consensus under all conditions
by slightly modifying Protocol I.
The modification will be made in step 3 of Protocol I, thereby
resulting in Protocol II.
Comparing Protocol I and Protocol II which follows, %it is easy to see that modifications are made in step 3.
the difference between the two only lies in the cases when an agent $i$ whose
gossip variable value at time $t$ equals that of its current preferred neighbor $i^*(t)$, at the
same time receives one or more requests to gossip. Under Protocol I,
agent $i$ gossips with that requesting neighbor whose label is closest to the front
of its neighbor queue at time $t$; the label $i^*(t)$ will still be in the front of the queue
at time $t+1$.
Under Protocol II, agent $i$ ignores all incoming requests to gossip at
time $t$ and moves the label $i^*(t)$ from the front to the end of the queue.

%The modified protocol is as follows.

\noindent{\bf Protocol II: }
Between clock times $t$ and $t+1$ each agent $i$ performs the steps enumerated below in the order
indicated. Although the agents' actions need not be precisely synchronized, it is understood that
for each $k\in\{1,2,3\}$ all agents complete step $k$ before any embark on step $k+1$.

\begin{enumerate}
\item Same as Protocol I
%{\bf 1st Transmission:}
%Agent $i$ sends its label $i$ and its gossip variable value $x_i(t)$ to its current
%preferred neighbor. At the same time agent $i$ receives the labels and corresponding gossip
%values from all of those neighbors which have agent $i$ as their current preferred neighbor.

\item Same as Protocol I
%{\bf 2nd Transmission:}
%Agent $i$ sends its current gossip value $x_i(t)$ to those neighbors which
%have agent $i$ as their current preferred neighbor.

\item {\bf Acceptances:}

\begin{description}
  \item[(a)] If $x_i(t)<x_{i^*(t)}(t)$ and agent $i$ has received at least one request
to gossip, then agent $i$ sends an acceptance to that particular requesting neighbor whose
label is closest to the front of the queue $\mathbf{q}_i(t)$.

  \item[(b)] If $x_i(t)\geq x_{i^*(t)}(t)$
  or agent $i$ has not received any request to gossip,
then agent $i$ does not send out an acceptance.
\end{description}

\item Same as Protocol I
%{\bf Gossip variable and queue updates:}
%\begin{description}
%\item[(a)] If agent $i$ sends an acceptance to or receives an acceptance from neighbor $j$,
%then agent $i$ gossips with neighbor $j$ by setting
%  $$x_i(t+1) = \frac{x_i(t)+x_j(t)}{2}$$
%  Agent $i$ updates its queue by moving $j$ from its current
%positions in $\mathbf{q}_i(t)$ to the end of the queue.

%  \item[(b)] If agent $i$ has not sent out an acceptance nor received one, then agent $i$ does not update
%the value of $x_i(t)$. In addition, $\mathbf{q}_i(t)$ is not updated except when agent $i$'s gossip
%value equals that of its current preferred neighbor. In this special case agent $i$
%moves the label $i^*(t)$ from the front to the end of the queue.
%\end{description}

\end{enumerate}
%{\bf Remark:}
%Comparing Protocol I and Protocol II, it is easy to see that modifications are made in step 3.
%The difference between the two protocols only lies in the cases when an agent $i$ whose
%gossip variable value at time $t$ equals that of its current preferred neighbor $i^*(t)$, at the
%same time receives one or more requests to gossip. Under Protocol I,
%agent $i$ gossips with that requesting neighbor whose label is closest to the front
%of its neighbor queue at time $t$; the label $i^*(t)$ will still be in the front of the queue
%at time $t+1$.
%Under Protocol II, agent $i$ ignores all incoming requests to gossip at
%time $t$ and moves the label $i^*(t)$ from the front to the end of the queue.

\noindent{\bf  Transmissions required:} At the end of step 1, agent $i$  knows
  $x_i(t)$ and $i^*(t)$ as well as the label and gossip  value of each neighbor
  which has   agent $i$ as its current preferred neighbor.
At the end of step 2, agent $i$ also knows  $x_{i^*(t)}(t)$. During step 1,
  each agent sends a transmission to its preferred neighbor so
 the total number of
 transmissions  required for all $n$ agents to complete step 1 is  $n$. During step 2,
  each neighbor of agent $i$ which has agent $i$ as its current preferred neighbor
 sends a
 transmission to agent $i$  so the total of
 transmissions  required for all $n$ agents to complete step 2 is  also $n$.
The total number of transmissions of all agents  required to complete step 3a is clearly
   no greater
   than $\frac{n}{2}$. Thus the total number of transmissions per iteration to carry out the protocol just described
     is no greater than $\frac{5}{2}n$. With a broadcasting protocol
  such as the one
   considered in \cite{alex}
   the total number of transmissions  per iteration  is
$$\sum_{i=1}^nn_i  = \sum_{i=1}^nd_i = n\left (\frac{1}{n}\sum_{i=1}^nd_i\right ) = nd_{{\rm avg}}$$
where $d_i$ is the degree of vertex $i$  and $d_{{\rm avg}}$ is the average vertex degree
$$d_{{\rm avg}}= \frac{1}{n}\sum_{i=1}^nd_i$$
 of the underlying  graph $\mathbb{A}$. Thus  for   gossip graphs with average vertex degree exceeding
  $\frac{5}{2}$,
 fewer transmissions are required  per iteration to do averaging with the protocol under consideration than are required
  per iteration to do averaging via broadcasting.

%Protocol II has the following properties:
It is possible to show that Protocol II is deadlock free.
To begin, let us note that at each time $t$, each label $i\in\{1,2,\ldots,n\}$
uniquely determines a sequence of labels
$$[i]_t=\{i_1,i_2,\ldots,i_{m(t)}\}$$
such that $i_1=i$, $i_{j+1}=i_j^*(t)$
for all $j\in\{1,2,\ldots,m(t)-1\}$, $i_1,i_2,\ldots,i_{m(t)-1}$ are distinct,
and $i_{m(t)}=i_k$ for some $k\in\{1,2,\ldots,m(t)-2\}$.
We call $[i]_t$ the sequence of {\em queue leaders} generated by $i$ at time $t$.
Note that $m(t)$ is a positive integer depending on time $t$ and always satisfies the inequalities
$2\leq m(t)-1\leq \delta$ where $\delta$ is the diameter of $\bbb{A}$.
We will sometimes simply write $[i]_t=\{i_1,i_2,\ldots,i_{m}\}$
for convenience with the understanding
that $m$ depends on time $t$.
The set of all possible sequences of queue leaders generated by $i$ is a finite set
because the number of agents in the group is finite.

\begin{lemma}
Suppose that all $n$ agents follow Protocol II. Then at each time $t$ either
$x_i(t)=x_{i^*(t)}(t)$ for some agent $i$, or a gossip must take place between two agents
whose gossip variables have different values.
\label{gossip}\end{lemma}

\noindent {\bf Proof of Lemma \ref{gossip}:}
Suppose that $x_i(t) \neq x_{i^*(t)}(t)$ for all $i\in\{1,2,\ldots,n\}$.
If a gossip takes place between agents $j$ and $k$, then either
$k=j^*$ or $j=k^*$; in either case $x_j(t)$ and $x_k(t)$ must have
different values because $x_i(t) \neq x_{i^*(t)}(t)$ for all $i\in\{1,2,\ldots,n\}$.
Thus to prove the lemma it is enough to show that gossip must take place
between two agents.

We claim that at least one agent must place a request to gossip.
To prove that this is so, suppose to the contrary. Then $x_i(t) < x_{i^*(t)}(t)$
for all $i\in\{1,2,\ldots,n\}$. In particular
$x_{i_1}(t)<x_{i_2}(t)<\ldots<x_{i_m}(t)$
where $\{i_1,i_2,\ldots,i_m\}=[i]_t$.
But $i_{m}=i_{m-1}^*(t)$ must equal some integer $i_j\in\{i_1,i_2,\ldots,i_{m-2}\}$
so $x_{i_{m-1}}(t)<x_{i_j}(t)$. But this is impossible because $j<m-1$.
Therefore at least one agent must place a request to gossip.

To complete the proof it is enough to show that among the agents who receive
requests to gossip at time $t$, at least one agent - say agent $k$ - does not
place a request to gossip. For if agent $k$ does not place a request, then
agent $k$ must gossip with that agent with label closest to the front of
$\mathbf{q}_k(t)$ who placed a request to gossip with agent $k$ at time $t$.
To prove that at least one agent receiving a gossip request at time $t$ does
not place a request to gossip at time $t$, assume the contrary.
Therefore suppose that every agent receiving a request to gossip at time $t$,
also places a request to gossip at time $t$. Let $i$ be the label of any agent
receiving a request to gossip at time $t$ and let $\{i_1,i_2,\ldots,i_m\}=[i]_t$.
Since agent $i_1=i$ and $i$ receives a request to gossip, it also must place a
request to gossip. Hence agent $i_2$ must receive a request to gossip.
Therefore agent $i_2$ must place a request to gossip at time $t$.
By this reasoning one concludes that all of the agents with label
$i_1,i_2,\ldots,i_{m-1}$ place requests to gossip at time $t$.
This implies that $x_{i_1}(t)>x_{i_2}(t)>\ldots>x_{i_m}(t)$.
But $i_{m}=i_{m-1}^*(t)$ must equal some integer $i_j\in\{i_1,i_2,\ldots,i_{m-2}\}$.
This means that $x_{i_{m-1}}(t)>x_{i_j}(t)$ which is impossible because $j<m-1$.
Therefore at least one agent which has received a request to gossip has not
placed a request to gossip.
\hfill
$\qed$

\begin{lemma}
Suppose that all $n$ agents follow Protocol II. Then a gossip must take place within
every $2d$ time steps,
where $d$ is the maximum vertex degree of $\bbb{A}$.
\label{dstep}\end{lemma}

\noindent {\bf Proof of Lemma \ref{dstep}:}
First, we claim that at least one agent must place a request to gossip within
every $d$ time steps. To prove that this is so, we need the following concept.
An agent $i\in\{1,2,\ldots,n\}$ is called {\em maximal} at time $t$ if its
gossip variable value is one of the largest at time $t$.
Let $\scr{M}(t)$ be the set of labels of all maximal agents at time $t$.
It is clear that $\scr{M}(t)=\{1,2,\ldots,n\}$
if and only if all $n$ agents have reached a consensus.
Thus for each time $t$, there must exist an agent $i\in\scr{M}(t)$
such that $i$ has at least one neighbor $j$ whose value $x_j(t)<x_i(t)$.
According to Protocol II, from time $t$ forward agent $i$ must repeatedly
either places a request to gossip with its preferred neighbor or updates its neighbor queue
by moving the label of its preferred neighbor from the front to the end of the queue until
it gossips with one of its neighbors.
For $i\in\{1,2,\ldots,n\}$, let $d_i$ be the
number of neighbors of agent $i$ or equivalently the degree of vertex $i$ in $\bbb{A}$.
In the worst case, when label $j$ is at the end of agent $i$'s neighbor queue at time $t$,
it will take at most $d_i\leq d$ time steps for agent $i$ to place a request to gossip.

To complete the proof it is enough to show that if at least one agent places a request to gossip
at time $\tau$, then a gossip must take place within the next $d$ time steps.
To prove that this is so, let $i$ be the label of any agent placing a request to
gossip at time $\tau$. By the same reasoning in the proof of Lemma \ref{gossip},
at least one agent - say agent $j$ - which receives a request to gossip does not
place a request to gossip at time $\tau$. Suppose $j$ receives a request to gossip
from its neighbor $k$.
From time $\tau$ forward, $k$ must repeatedly
places a request to gossip with $j$ until $j$ gossips with one of its neighbors.
By the protocol, $j$ does not gossip if and only if $k$ is not $j$'s preferred neighbor
and $j$'s gossip value equals that of its preferred neighbor;
if it is so, $j$ updates its neighbor queue and the label $k$ advances closer to
the front of the queue.
In the worst case, when label $k$ is at the end of $j$'s neighbor queue at time $\tau$,
it will take at most $d_j\leq d$ time steps for label $k$ to advance to the front of
$j$'s queue. This means that agents $j$ and $k$ are guaranteed to gossip at least
once within $d$ time steps after time $\tau$.
\hfill
$\qed$

It is also possible to show that every sequence of gossip vectors generated by
Protocol II converges to the desired limit point exponentially fast.

\begin{theorem}
Suppose that all $n$ agents adhere to Protocol II. Then there is a finite time $T$,
not depending on the values of gossip variables, such that every sequence of gossip
vectors $x(t)$, $t>0$ generated by Protocol II is repetitively complete with period
no greater than $T$.
\label{complete2}\end{theorem}

To prove Theorem \ref{complete2}, we need a few ideas.
First note that
step 4 of the protocol stipulates that agent $i$ must update its queue
whenever its current gossip value equals that of
its current preferred neighbor. We say that agent $i$ {\em gossips virtually}
with neighbor $j$ at time $t$ if $i^*(t)=j$ and the current gossip
values of both agents are the same.
%It is easy to see that virtual gossiping is not a symmetric relation.
It is worth noting that when agent $i$ gossips virtually with neighbor $j$,
$j$ may not gossip virtually with $i$.
Also note that each agent can gossip virtually with at most one neighbor at one clock time.
If an agent  gossips virtually with its current preferred neighbor,
it does not gossip with any other neighbor.
Thus each agent can gossip or virtually gossip with at most one neighbor
at one clock time.
If agent $i$ gossips or gossips virtually with neighbor $j$ at time $t$,
then agent $i$ updates its neighbor queue by moving the label $j$ from its current position in
$\mathbf{q}_i(t)$ to the end of the queue.

To proceed, we say that an agent has completed a {\em round} of gossiping after
it has gossiped or virtually gossiped  with
 each neighbor in $\scr{N}_i$ at least once.   Thus the finite
 sequence of primitive gossiping matrices corresponding to a finite sequence
 of multi-gossips and virtual multi-gossips
for the entire group  which has occurred over an interval of length $T$,
  will be  complete
 if over the same period  each agent in the group completes a round.
 Thus Theorem \ref{complete2} will be true if
 every agent completes a round in a number of iterations no larger than $T$.
The following proposition asserts that this is in fact the case.

\begin{proposition}
Suppose that all $n$ agents adhere to Protocol II. Then there is a finite time $T$,
not depending on the values of gossip variables, such that within $T$ iterations every
agent will have gossiped or virtually gossiped at least once  with each of its neighbor.
\label{reed}\end{proposition}
To prove this proposition we will make use of the following two lemmas.

%To prove Theorem \ref{complete2}, we need several preliminary results.
%First, it is possible to show that for each $i$ and each time $t$ at which $x_i(t) \neq x_{i^*(t)}(t)$ there is a
%finite time $T$, not depending on the initial values of the gossip variables, within which agents $i$ and
%$i^*(t)$ will gossip. %To this end we need several preliminary results.
%This and the manner in which the queues are updated can then be used to show
%that any gossip sequence generated by this protocol must be repetitively complete with a period
%bounded above by a finite number not depending on initial values of the gossip variables.

%\begin{lemma} Suppose that all $n$ agents adhere to Protocol II. Then for
%each $i\in\{1,2,\ldots, n\}$  there is a finite time $T_i$  with the following property.
% For each $t\geq 0$ at which $x_i(t) \neq x_{i^*(t)}(t)$ there is a
%  time  $\tau\in [t,t+T_i]$ at which either $x_i(\tau) = x_{i^*(t)}(\tau)$
%or agents $i$ and $i^*(t)$ gossip.
%\label{pizza2}\end{lemma}
%To prove this lemma we need several preliminary results.

%For each $i\in\{1,2,\ldots,n\}$  and each $t\geq 0$
%there always exits a unique path in $\bbb{A}$,
%$(i,i_1),(i_1,i_2),\ldots$, $(i_{m-1},i_m)$, such that $i^*(t)=i_1$,
%$i_j^*(t)=i_{j+1}$, $j\in\{1,2,\ldots,m-1\}$, and $i_m^*(t)\in\{i,i_1,i_2,\ldots,i_{m-1}\}$.

\begin{lemma}
Suppose that all $n$ agents adhere to Protocol II.
Then for each
$i\in\{1,2,\ldots,n\}$ there is a finite time $T$ with the following property:
For each $t\geq 0$ let $[i]_t=\{i_1,i_2,\ldots,i_m\}$.
There must be a time  $\tau\in [t,t+T]$ and one integer $j\in\{1,2,\ldots,m-1\}$
such that agent $i_j$ either gossips or virtually gossips with agent $i_{j+1}$.
\label{go1}\end{lemma}

\noindent {\bf Proof of Lemma \ref{go1}:}
%Given any $i\in\{1,2,\ldots,n\}$ and $t\geq 0$.
%Let $[i]_t=\{i_1,i_2,\ldots,i_m,i_{m+1}\}$.
For any two successive times $t$ and $t+1$,
$[i]_{t+1} = [i]_t$ if and only if there is no integer $j\in\{1,2,\ldots,m-1\}$
such that agent $i_j$ either gossips or virtually gossips with agent $i_{j+1}$ at time $t$.
In other words, $[i]_{t+1} \neq [i]_t$ if and only if there is at least
one integer $j\in\{1,2,\ldots,m-1\}$
such that agent $i_j$ either gossips or virtually gossips with agent $i_{j+1}$ at time $t$.
Thus to prove the lemma, it is enough to show that for each
$i\in\{1,2,\ldots,n\}$ there is a finite time $T$ with the following property:
For each $t\geq 0$
there is must be a time  $\tau\in [t,t+T]$ such that
$[i]_{\tau+1} \neq [i]_\tau$.

%\begin{lemma}
%Suppose that all $n$ agents adhere to Protocol IV. Then for each
%$i\in\{1,2,\ldots,n\}$ there is a finite time $T$ with the following property.
%For each $t\geq 0$ at which $[i]_t=\{i_1,i_2,\ldots,i_m,i_{m+1}\}$
%there is a time  $\tau\in [t,t+T]$ at which
%either $x_{i_j}(\tau) = x_{i_j^*(t)}(\tau)$ or
%agents $i_j$ and $i_j^*(t)$ gossip for some $j\in\{1,2,\ldots,m\}$.
%\label{go1}\end{lemma}

%\begin{lemma}
%Suppose that all $n$ agents adhere to Protocol II. Then for each
%$i\in\{1,2,\ldots,n\}$ there is a finite time $T$ with the following property.
%For each $t\geq 0$
%there is a time  $\tau\in [t,t+T]$ such that
%$[i]_{\tau+1} \neq [i]_\tau$.
%\label{go1}\end{lemma}

%Let $[i]_t=\{i_1,i_2,\ldots,i_m,i_{m+1}\}$.
Since $i_{m}=i_k$ for some $k\in\{1,2,\ldots, m-2\}$,
 $[i_k]_t=\{i_k,i_{k+1},\ldots,i_{m-1},i_k\}$.
%If there is no $j\in\{1,2,\ldots,m\}$ such that ......
%$[i]_{t+1}=[i]_k$ and $[i_k]_{t+1}=[i_k]_t$.
%Then $(i_k,i_{k+1}),(i_{k+1},i_{k+2}),\ldots,(i_{m-1},i_m)$
%is a path in $\bbb{A}$ such that $i_j^*(t)=i_{j+1}$, $j\in\{k,k+1,\ldots,m-1\}$,
%and $i_m^*(t)=i_k$.
Suppose that $x_{i_j}(t)\neq x_{i_j^*(t)}(t)$ for all $j\in\{k,k+1,\ldots,m-1\}$.
Using the same arguments in the proof of Lemma \ref{gossip},
at least one agent in $\{i_k,i_{k+1},\ldots,i_{m-1}\}$ receiving a request to
gossip at time $t$ does not place a request to gossip at time $t$.
In other words, there is at least one integer $j\in\{k,k+1,\ldots,m-1\}$ such that
$i_j$ places a request to gossip with $i_j^*(t)$ at time $t$
and $i_j^*(t)$ does not place a request to gossip at time $t$.

Note that if $[i_k]_{t+1} \neq [i_k]_t$ then $[i]_{t+1} \neq [i]_t$.
Thus to complete the proof it is enough to show that
there is a finite time $T$ such that
there is a time  $\tau\in [t,t+T]$ at which $[i_k]_{\tau+1} \neq [i_k]_\tau$.
At time $t$, suppose $i_j\in\{i_k,i_{k+1},\ldots,i_{m-1}\}$
places a request to gossip with $i_j^*(t)$ and $i_j^*(t)$ does not place
a request to gossip.
$i_j^*(t)$ does not gossip with
$i_j$ if and only if $i_j^*(t)$ at the same time receives a request to gossip
from some other
neighbors whose labels are closer to the front of $\mathbf{q}_{i_j^*(t)}$.
In this case $i_j^*(t)$ gossips with that requesting
neighbor whose label is closest to the front of $\mathbf{q}_{i_j^*(t)}$;
the label of $i_{j}$ thus advances closer to the front of $\mathbf{q}_{i_j^*(t)}$.
Therefore at time $t$ either $[i_k]_{t+1} \neq [i_k]_t$, or there is at least
one integer $j\in\{k,k+1,\ldots,m-1\}$ such that the label of $i_{j}$ advances
closer to the front of $\mathbf{q}_{i_j^*(t)}$.
If the latter situation occurs, $[i_k]_{t+1} = [i_k]_t$.
Then by the same reasoning as before,
at time $t+1$ either $[i_k]_{t+2} \neq [i_k]_{t+1}$, or there is at least
one integer $j\in\{k,k+1,\ldots,m-1\}$ such that the label of $i_{j}$ advances
closer to the front of $\mathbf{q}_{i_j^*(t+1)}$.
Let $d_{i}$ denote the number of neighbors of agent $i$ or equivalently the degree of
vertex $i$ in $\bbb{A}$.
In the worst case, when label $i_j$ is at the end of $\mathbf{q}_{i_j^*(t)}$
at time $t$ for all $j\in\{k,k+1,\ldots,m-1\}$,
it will take at most
$\sum_{j=k}^{m} d_{i_j}$ successive times for
each label $i_j\in\{i_k,i_{k+1},\ldots,i_{m-1}\}$ to advance to the front of $\mathbf{q}_{i_j^*(t)}$.

If at some time $\tau>t$, $[i_k]_\tau=[i_k]_t=\{i_k,i_{k+1},\ldots,i_{m-1},i_k\}$ and
each label $i_j\in\{i_k,i_{k+1},\ldots,i_{m-1}\}$ is at the front of $\mathbf{q}_{i_j^*(t)}$,
there is at least one integer $j\in\{k,k+1,\ldots,m-1\}$ such that
either $x_{i_j}(\tau) = x_{i_j^*(t)}(\tau)$ or
agents $i_j$ and $i_j^*(t)$ gossip at time $\tau$.
In other words, this implies that $[i_k]_{\tau+1} \neq [i_k]_\tau$.
Therefore for each $i\in\{1,2,\ldots,n\}$ and $t\geq 0$, a time $\tau\geq t$
at which $[i]_{\tau+1}\neq [i]_\tau$ is guaranteed
in any time interval containing
$\sum_{j=k}^{m} d_{i_j}$ times.
\hfill
$\qed$

%For those $t$ such that $[i]_{t+1}\neq [i]_t$,
%if $[i]_t$ repeats itself at some time after $t$,
% more can be said.

\begin{lemma}
Suppose that all $n$ agents adhere to Protocol II.
For each $t\geq 0$ and $i\in\{1,2,\ldots,n\}$,
let $t+T$ be the first time that $[i]_t$ repeats itself.
If $T>1$, then there is a
  time  $\tau\in [t,t+T]$ at which agents $i$ and $i^*(t)$ either
gossip or gossip virtually.
\label{go}\end{lemma}

\noindent {\bf Proof of Lemma \ref{go}:}
Let $[i]_t=\{i_1,i_2,\ldots,i_m,i_{m+1}\}$ in which $i_1=i$ and
$i_{j+1}=i_j^*(t)$ for all $j\in\{1,2,\ldots,m\}$.
With $T>1$, that $t+T$ is the first time that $[i]_t$ repeats itself
implies that there is a time $t_1\in[t,t+T-1]$ at which $[i]_{t_1+1}\neq [i]_{t_1}$.
This means that there is an integer $j\in\{1,2,\ldots,m\}$
such that either $x_{i_j}(t_1) = x_{i_{j+1}}(t_1)$ or agents $i_j$ and $i_{j+1}$ gossip at time $t_1$.
If $j=1$ then the lemma is clearly true. Now suppose that $j>1$.
The label of $i_{j+1}$ is moved to the end of the neighbor queue of $i_j$ at time $t_1$.
Since $[i]_{t+T}=[i]_t$, %$i_j^*(t+T)=i_j^*(t)=i_{j+1}$.
the label of $i_{j+1}$ advances to the front of the neighbor queue of $i_j$ at time $t+T$.
This implies that $i_j$ has completed a round of of gossiping
during the time interval $[t_1,t+T]$. Since $i_{j-1}$ is a neighbor of $i_j$,
there is a time  $\tau\in [t_1,t+T]$ at which either $x_{i_{j-1}}(\tau)=x_{i_j}(\tau)$
or agents $i_{j-1}$ and $i_j$  gossip. %Since $i_j=i_{j-1}^*(t)$, it is
%equivalent to that there is a time  $\tau\in [t,t+T]$ at which either
%$x_{i_{j-1}}(\tau) = x_{i_{j-1}^*(t)}(\tau)$
%or agents $i_{j-1}$ and $i_{j-1}^*(t)$ gossip.
A simple induction
thus proves that there is a
  time  $\tau\in [t,t+T]$ at which either $x_{i_1}(\tau) = x_{i_2}(\tau)$
or agents $i_1$ and $i_2$ gossip.
\hfill
$\qed$

%\begin{lemma} Suppose that all agents
%in the group adhere to Protocol IV.
%Then for each $i$ and each time $t$ at which $x_i(t)\neq x_{i^*(t)}(t)$
%there is a finite time $T$, not depending on the initial values of the gossip variables,
%within which either $x_i=x_{i^*}$ or agent $i$ and $i^*(t)$ will gossip.
%\label{bound}\end{lemma}

\noindent
We are now in a position to prove Proposition \ref{reed}.

%\noindent {\bf Proof of Lemma \ref{pizza2}:}
%Let $[i]_t=\{i_1,i_1,\ldots,i_m,i_{m+1}\}$.
%In view of Lemma \ref{go1}, there is a finite time $T$ with the following property.
% For each $t\geq 0$
%there is a time  $\tau\in [t,t+T]$ such that
%$[i]_{\tau+1} \neq [i]_\tau$.
%In other words, $[i]$ is guaranteed to change at least once within any interval
%containing $T$ times.

%Recall that the set of all possible sequences of queue leaders generated by $i$ is a finite set.
%Let $N$ denote the number of all possible sequences of queue leaders generated by $i$.
%Then within a finite time $T_i=(N+1)T$, there exist at least $N+1$ distinct times
%$\tau_1,\tau_2,\ldots,\tau_{N+1}\in[t,t+T_i]$
%such that $[i]_{\tau_{k}+1}\neq [i]_{\tau_{k}}$ for all $k\in\{1,2,\ldots,N+1\}$.
%Among these $N+1$ distinct times, there must exist
%two times $t_1,t_2\in\{\tau_1,\tau_2,\ldots,\tau_{N+1}\}$
%such that $[i]_{t_1}=[i]_{t_2}$ and $t_2-t_1>1$.
%In view of Lemma \ref{go}, this ensures that there is a
%  time  $\tau\in [t,t+T_i]$ at which either $x_i(\tau) = x_{i^*(t)}(\tau)$
%or agents $i$ and $i^*(t)$ gossip.
%$\qed$

\noindent
{\bf Proof of Proposition \ref{reed}:}
For each $t\geq 0$ and $i\in\{1,2,\ldots,n\}$,
let $[i]_t=\{i_1,i_1,\ldots,i_m,i_{m+1}\}$.
In view of Lemma \ref{go1}, there is a finite time $T$,
not depending on the initial values of gossip variables, with the following property:
There must be a time  $\tau\in [t,t+T]$ such that
$[i]_{\tau+1} \neq [i]_\tau$.
In other words, $[i]_t$ is guaranteed to change at least once within any interval
containing $T$ times.

Recall that the set of all possible sequences of queue leaders generated by $i$ is a finite set.
Let $N$ denote the number of all possible sequences of queue leaders generated by $i$.
Then within a finite time $T_i=(N+1)T$, there exist at least $N+1$ distinct times
$\tau_1,\tau_2,\ldots,\tau_{N+1}\in[t,t+T_i]$
such that $[i]_{\tau_{k}+1}\neq [i]_{\tau_{k}}$ for all $k\in\{1,2,\ldots,N+1\}$.
Among these $N+1$ distinct times, there must exist
two times $t_1,t_2\in\{\tau_1,\tau_2,\ldots,\tau_{N+1}\}$
such that $[i]_{t_1}=[i]_{t_2}$ and $t_2-t_1>1$.
In view of Lemma \ref{go}, this ensures that there is a
  time  $\tau\in [t,t+T_i]$ at which agents $i$ and $i^*(t)$ either gossip
  or gossip virtually.
  Once agent $i$ gossips or virtually gossips with $i^*(t)$, the protocol
  stipulates that agent $i$ must move the label $i^*(t)$ to the end of its queue.
Thus agent $i$ is guaranteed to complete a round of gossiping
within $d_{i}T_i$ time steps, where $d_i$ is the number of neighbors of agent $i$.
A finite gossiping sequence for the entire group which has occurred over an interval of
length $T$ will be complete if over the same period each agent in the group completes a round.
Therefore the length of the time interval large enough for all agents in the group to complete a round of
gossiping is the maximum of the times $d_iT_i$, $i\in\{1,2,\ldots,n\}$.
%We are led to the following result.
\hfill
$\qed$

%For each agent $i$ and time $t$,
%let $i_1=i$ and $\bbb{P}(t)=(i_1,i_2),(i_2,i_3),\ldots,(i_{m-1},i_m)$.
%In view of the proof of Lemma \ref{bound}, $T_i\leq N_i\tau+1$, where $N_i$ denotes
%the number of all possible waiting paths generated by $i$ with starting edge
%$(i_1,i_2)$ and $\tau$ is a finite time within which at least one gossip or pseudo gossip
%takes place between $i_j$ and $i_j^*(t)$ for some $j\in\{1,2,\ldots,m\}$.
%$\tau\leq D\Delta$ where $D$ and $\Delta$ are the
%diameter and degree of $\bbb{A}$ respectively.
%Note that $(i_1,i_2)$ is a possible waiting path for which $\tau< D\Delta$.
%Thus $T_i\leq N_i\tau$.
%The length of each possible waiting path cannot exceed the diameter of $\bbb{A}$.
%Thus $m\leq D+1$.
%Consider a sequence of $D+2$ agents, $i_1,i_2,\ldots,i_{D+2}$.
%If $i_1$, $i_2$ are fixed and $i_3,\ldots,i_{D+2}$ can be any integer
%$j\in\{1,2,\ldots,n\}$, then there are totally $n^D$ such sequences.
%Clearly $i_{D+2}\in\{i_1,i_2,\ldots,i_{D+1}\}$. Thus for each possible
%waiting path $\bbb{P}(t)=(i_1,i_2),(i_2,i_3),\ldots,(i_{m-1},i_m)$, there always
%exists a sequence of $D+2$ agents whose first $m+1$ agents are
%$i_1,i_2,\ldots,i_m,i_m^*(t)$. This implies that $N_i\leq n^D$.
%Therefore $\max_i (d_iT_i)\leq \max_i (d_i\Delta D n^D) = \Delta^2 D n^D$.

In the light of Theorem \ref{complete2} it is of interest to derive an explicit expression of $T$.
But the tight bound of $T$ has so far
eluded us except for the special case when $\bbb{A}$ is a tree.

\begin{theorem}
Suppose that the graph of allowable gossips $\bbb{A}$ is a tree with $n>1$ agents
and that all $n$ agents adhere to Protocol II. Then every sequence of gossip vectors
$x(t)$, $t>0$ generated by Protocol II is repetitively complete with period no greater
than the number of edges of $\bbb{A}$.
\label{tree}\end{theorem}

To prove Theorem \ref{tree}, we need to generalize slightly the idea of a round.
%For any nonempty subset $\scr{M}\subset\scr{N}_i$, we say that
%agent $i$ has completed a {\em round of gossiping with $\scr{M}$} within $T$ clock times if
%over the same period and for each $j\in\scr{M}$, agent $i$ gossips or virtually gossips with
%agent $j$ at least once.
For any nonempty subset $\scr{M}\subset\scr{N}_i$, we say that
agent $i$ has completed a {\em round of gossiping with $\scr{M}$} if
agent $i$ has gossiped or virtually gossiped with each neighbor in $\scr{M}$ at least once.

\noindent {\bf Proof of Theorem \ref{tree}:}
A finite gossiping sequence for the entire group which has occurred over an interval of
length $T$ will be complete if over the same period each agent in the group completes a round
of gossiping.
Thus to prove Theorem \ref{tree} it is sufficient to
establish the following claim.

\noindent{\bf Claim:}
Let $\bbb{T}$ be any fixed tree with $n>1$ agents.
For each time $t\geq 0$ and each $i\in\{1,2,\ldots,n\}$,
agent $i$ will complete a round
of gossiping within $n-1$ iterations starting at time $t$.

\noindent This claim will be proved by induction on $n$.

Suppose $n = 2$ in which case $\bbb{T}$ has one edge.
Let $i$ and $j$ be the two agents. It is clear that $i^*(t)=j$ and
$j^*(t)=i$ for all $t>0$. If $x_i(t)\neq x_j(t)$ and, without loss of generality,
assume that $x_i(t)>x_j(t)$, then only $i$ places a request to gossip with $j$
and $j$ must accepts the request.
Thus in this case the claim is true.

Now suppose that the claim holds for all $n$ in the range $2\leq n\leq m$ where $m$ is
a positive integer greater than 1. Let $\bbb{T}$ be any fixed tree with $m + 1$ agents.

Let $i\in\{1,2,\ldots,n\}$ be any agent in $\bbb{T}$ and
let $d$ be the number of neighbors of agent $i$
or equivalently the degree of vertex $i$ in $\bbb{T}$.
Let $v_1,v_2,\ldots,v_d$ denote all $d$ neighbors of agent $i$.
It is possible to decompose $\bbb{T}$ into a set of $d$ spanning subgraphs
which are linked together at vertex $i$.\footnote{
A finite set of spanning subgraphs $\bbb{G}_1,\bbb{G}_2,\ldots,\bbb{G}_m$ of
a simple graph $\bbb{G}$ is {\em linked together at vertex $v$} if $v$
has degree 1 in each subgraph and, in addition, if $v$ is the only ``non-isolated''
vertex of each subgraph which every possible pair of subgraphs have in common.
By a {\em non-isolated} vertex of a simple graph $\bbb{H}$ is meant any vertex
in $\bbb{H}$ with positive degree.}
Towards this end let $\scr{E}_i$ denote the set of $d$ edges of $\bbb{T}$
which are incident on vertex $i$, and label them 1 through $d$.
For each $j\in\{1,2,\ldots,d\}$ let $\bbb{G}_j$ denote the spanning subgraph of $\bbb{T}$ which
results when edge $j$ is deleted from $\bbb{T}$.
Let $\bbb{T}_j$ be that connected component of $\bbb{G}_j$ which contains vertex $v_i$.
%{\color{red} contains vertex $v_i$}
Clearly each $\bbb{T}_j$
is a tree. Let $n_j$ denote the number of agents of $\bbb{T}_j$.
Then $2\leq n_j\leq m$ for all $j\in\{1,2,\ldots,d\}$ and
$\sum_{j=1}^d n_j=m$.

Without loss of generality, suppose that $\mathbf{q}_i(t)=\{v_1,v_2,\ldots,v_d\}$.
If $v_1^*(t)=i$, then either $x_i(t)=x_{v_1}(t)$ or agents $i$ and $v_1$ gossip at time $t$.
If $v_1^*(t)\neq i$ and assume, in the worst case, that the label $i$ is
at the end of $\mathbf{q}_{v_1}(t)$, then
$i$ will advance to the front of the queue of $v_1$
after $v_1$ has completed a round of gossiping with $\scr{N}_i\setminus\{i\}$.
Note that in $\bbb{T}_1$ the neighbor set of $v_1$ consists of all the neighbors of
$v_1$ in $\bbb{T}$ except for $i$.
Thus in the worst case the label $i$ will advance to the front of the queue
of agent $v_1$ after $v_1$ has completed a round of gossiping in $\bbb{T}_1$.
Since $\bbb{T}_1$ is a tree whose number of agents is in the range $2\leq n_1\leq m$,
$v_1$ will complete a round of gossiping in at most $n_1-1$ clock times by the inductive hypothesis.
Therefore over a period of at most $n_1$ clock times there is at least one time
$\tau$ such that either $x_i(\tau)=x_{v_1}(\tau)$ or agents $i$ and $v_1$ gossip at time $\tau$.
Then the label of agent $v_2$ advances to the front of the queue of agent $i$.
By the same reasoning as before, it will take at worst, an additional $n_2$ successive
clock times for the label of agent $v_3$ to advance to the front of the queue of agent $i$
\{i.e., there is a time $\bar{\tau}$ such that either $x_i(\bar{\tau})=x_{v_2}(\bar{\tau})$
or agents $i$ and $v_2$ gossip at time $\bar{\tau}$\}.
In other words, agent $i$ is guaranteed to have completed a round of gossiping with
$\{v_1,v_2\}$ at least once within any time interval containing no more than $n_1+n_2$ clock times.
By repeating this argument for all labels in the queue $\mathbf{q}_i(t)$,
one reaches the conclusion that agent $i$ is guaranteed to complete a round of gossiping
within any time interval containing at most $\sum_{j=1}^d n_j=m$ clock times.
By induction, the claim is established and the proof is complete.
\hfill
$\qed$

\subsection{An Accelerated Protocol}\label{acc}

An important rule of gossiping is that during a gossiping process
each agent is allowed to gossip with at most one of its neighbors at one clock time.
There is no such restriction on virtual gossips.
Thus to improve the convergence rate of the protocol in the preceding section,
a natural idea is to let each agent
gossip virtually with as many as neighbors as possible at the same time.
%Here the idea of a virtual gossip is generalized. We say that agent $i$ {\em gossip virtually}
%with neighbor $j$ at time $t$ if the current gossip values of both agents are the same.
%%Such an accelerated protocol is as follows and has been reported in \cite{cdc11req}.

\noindent{\bf Protocol III: }
Between clock times $t$ and $t+1$ each agent $i$ performs  the steps enumerated below  in the order indicated.
 Although the agents'  actions need not be
precisely synchronized, it is understood that  for each
 $k\in\{1,2,3\}$  all agents complete step $k$ before any embark on step $k+1$.

\begin{enumerate}
\item Same as Protocol I
%{\bf  1st Transmission:}
%Agent $i$  sends
% its label $i$ and its gossip    value $x_i(t)$
% to  its current preferred neighbor.  At the same time agent $i$ receives
%  the  labels and corresponding gossip values from   all of
%   those neighbors which have agent $i$ as their current preferred neighbor.

\item Same as Protocol I
%{\bf  2nd Transmission:}
%Agent  $i$  sends its current gossip value
%$x_i(t)$
%to those neighbors which have agent $i$ as their current preferred neighbor.

\item Same as Protocol II
%{\bf Acceptances:} \begin{enumerate}
%\item If agent $i$ has not placed a request to gossip
% but has received at least one request to gossip,  then agent $i$ sends an
% acceptance to  that particular  requesting neighbor
% whose label is closest to the
%  front  of the
% queue $\mathbf{q}_i(t)$.

% \item If agent $i$ has either  placed a request to gossip   or
%  has not received any requests to gossip,
%   then agent $i$ does not send
%  out an acceptance. \end{enumerate}

%\item {\bf Acceptances:} \begin{enumerate} \item If $x_i(t)\leq x_{i^*(t)}(t)$ and agent $i$ has at
%least one dominant requester, then
% agent $i$ sends
%an acceptance to that particular dominant requester whose label is closest to the front
% of $q_i(t)$. \item  If $x_i(t)> x_{i^*(t)}(t)$ or if agent $i$ has no   dominant requester,
 %agent $i$ does not send out an
 %acceptance. \end{enumerate}

%\item If agent $i$ sends  an acceptance to co-preference $j$, then co-preference $j$ sends $x_j$ to agent $i$.
\item  {\bf Gossip variable and queue updates:}\begin{enumerate}

\item If agent $i$ either sends an acceptance to or receives an acceptance from  neighbor $j$,
then agent $i$ gossips with neighbor $j$ by setting
 $$x_i(t+1) = \frac{x_i(t)+x_{j}(t)}{2}$$%\hspace{.5in} {\rm and}\hspace{.5in}  x_j(t+1) = \frac{x_i(t)+x_{j}(t)}{2}$$
%respectively.
Agent $i$ updates its queue  by moving  $j$ and the labels of all of its current receivers $k$,
if any,
for which $x_k(t) = x_i(t)$
from their current positions
in
  $\mathbf{q}_i(t)$ to the end of the queue  while maintaining their relative order.

%\item If agent $i$ sends an  acceptance to neighbor  $j$,  then agent $i$ sets
 %$x_i(t+1) = \frac{x_i(t)+x_{j}(t)}{2}$ and  updates $q_i(t)$ by moving
 % $j$ from its current position in  $q_i(t)$ to the end of the queue

 %\item  If agent
%$i$  receives  an acceptance from its preferred neighbor  $j$,
% then agent $i$ sets
% $x_i(t+1) = \frac{x_i(t)+x_{j}(t)}{2}$ and  updates $q_i(t)$ by moving
 % $j$ from its  position in front of   $q_i(t)$ to the end of the queue.

  \item If agent $i$ has not sent out an acceptance nor received one, then agent $i$ does not
  update the  value of $x_i(t)$.
In addition,  $\mathbf{q}_i(t)$ is not updated except
    when agent $i$'s gossip  value equals that of
    at least one of its current receivers. In this special case agent $i$ moves the
    labels of all
    of its current receivers
    $k$ for which $x_k(t) = x_i(t)$  from their current positions in
  $\mathbf{q}_i(t)$ to the end of the queue, while maintaining their relative order.

\end{enumerate}

\end{enumerate}

%In summary,
%\begin{itemize}
%\item  For agent $i$ to place a request to gossip,  the current value of its gossip variable must be larger
%than that of its current preferred neighbor.
% \item  For a gossip to occur   between two agents $i$ and $j$ at time $t$,
%one -- say $i$ -- must be the current  preferred
%neighbor  of the other \{i.e., $i=j^*(t)$\},
%$x_j(t)$ must be larger than $x_{i}(t)$, and $j$ must be the label of the
%neighbor of agent $i$ with highest priority which is placing a request
% to gossip with agent $i$.
% \item For agent $i$ to update its queue it must either gossip with a neighbor $j$
%  or, if not,   it's current gossip
% value must equal that of at least one of its receivers.
%\end{itemize}

Protocol III is expected to solve the distributed averaging problem faster than
Protocol II since Protocol III allows agents to ``gossip virtually'' with more than one
neighbors at one time while  Protocol II dose not.
Faster convergence of Protocol III was illustrated in \cite{cdc12gossip} by simulation
(see Section V in \cite{cdc12gossip}).

It is also possible to derive a tight bound on the convergence rate of Protocol III
for general allowable gossip graphs.

\begin{theorem}
Suppose that all $n$ agents follow Protocol III.
Then for any connected allowable gossip graph $\bbb{A}$,
every sequence of gossip vectors  $x(t)$, $t>0$  generated by
 Protocol III is  repetitively complete with period no greater than the number
 of edges of $\bbb{A}$.
\label{ji}\end{theorem}
%A proof of this theorem can be found in \cite{cdc11req}.

To prove Theorem \ref{ji}, we need to generalize slightly a few  ideas. First note that
step 4 of the protocol stipulates that agent $i$ must update its queue whenever its current gossip value equals that of
on of its neighbors. We say that agent $i$ {\em gossips virtually} with neighbor $j$ at time $t$ if the current gossip
values of both agents are the same. Note that while an agent can gossip with at most one agent at time $t$, it can
gossip virtually with as many as $n_i$ at the same time.
%To proceed,  we need to generalize slightly the idea of a round.
We say that an agent has completed a round of gossiping after it has gossiped or virtually gossiped  with
 each neighbor in $\scr{N}_i$ at least once.   Thus the finite
 sequence of primitive gossiping matrices corresponding to a finite sequence of multi-gossips and virtual multi-gossips
for the entire group  which has occurred over an interval of length $T$,
  will be  complete
 if over the same period  each agent in the group completes a round. Thus Theorem \ref{ji} will be true if
 every agent completes a round in a number of iterations no larger than the number of edges of $\bbb{A}$.
The following proposition asserts that this is in fact the case.

\begin{proposition} Let $m$ be the number of edges in $\bbb{A}$.  Then within $m$ iterations every
agent will have gossiped or virtually gossiped at least once  with each of its neighbor.
\label{reeds}\end{proposition}
To prove this proposition we will make use of the following two lemmas.

\begin{lemma}
Suppose that all $n$ agents follow Protocol III. Then at each time $t$, at least one gossip or
 virtual gossip must occur.
\label{pizza}\end{lemma}

\begin{lemma} Let $t$ be fixed and suppose that $\bbb{G}$ is
 a spanning subgraph of $\bbb{A}$ with at least one edge.
  For each $i\in\{1,2,\ldots,n\} $ write $\scr{N}_i$
 for the set of labels of the vertices
adjacent to vertex $i$ in $\bbb{A}$ and $\scr{M}_i$ for the set of
labels of the vertices
adjacent to vertex $i$ in $\bbb{G}$.
    Let  $\scr{N}_i-\scr{M}_i$ denote the
   complement of $\scr{M}_i$ in $\scr{N}_i$.
    Suppose that for each
    $i\in\{1,2,\ldots, n\}$, each label in $\scr{M}_i$,
if any, is  closer to the front of $\mathbf{q}_i(t)$ than are all the labels
in $\scr{N}_i-\scr{M}_i$. Then
 there must be an edge $(i,j)$ within $\bbb{G}$  such that at time $t$,
 neighboring agents $i$ and $j$ either gossip or
gossip virtually.\label{lib} \end{lemma}

\noindent
We will prove lemma \ref{lib} first.
%In its proof,  we use the following notation.  For $i\in\{1,2,\ldots,n\}$,
%  $[i]_t$ denote the sequence of $w>1$ distinct integers
%$[i]_t = \{i_1,i_2,\ldots,i_{w}\}$ satisfying $i_1=i$, $i_{j+1} = i_j^*(t),\;j\in\{1,2,
%\ldots,w-1\}$ and $i_{w}^*(t) = i_k$ for some $k\in\{1,2,\ldots,w-1\}$.
%Note that $w$ depends on time and always satisfies the inequalities
% $2\leq w\leq \delta $ where $\delta $ is the diameter of $\bbb{A}$.

%Questions:  Is a graph with one vertex connected?

\noindent{\bf Proof of Lemma \ref{lib}:}  Let  $\scr{J}$ denote the set of labels of
 all agents $i$  for which  $\scr{M}_i$ is nonempty.
Since $\bbb{G}$ has at least one edge, $\scr{J}$ is nonempty.
 Fix  $i\in\scr{J}$. We claim that
 $i^*(t)$ must be in $\scr{M}_i$. If it were not, it would have
  to be further back in $\mathbf{q}_i(t)$ than  the labels in $\scr{M}_i$
  and this would contradict
   the fact that $i^*(t)$ is in the front of $\mathbf{q}_i(t)$. Therefore $i^*(t)\in \scr{M}_i$. This implies that
    $(i,i^*(t))$ is an edge in
   $\bbb{G}$.  Hence   $\scr{M}_{i^*(t)}$ must be nonempty so  $i^*(t)$
   must also be in $\scr{J}$. From this it follows that for each $i\in\scr{J}$, all of the labels in
   $[i]_t$ are also  in $\scr{J}$.

To proceed,  suppose that
 $x_i(t) =  x_{i^*(t)}(t)$
   for some $i\in\scr{J}$. Then agent $i$ has not placed a request. If agent $i$ receives a request,
   then agent $i$ must send an acceptance because of 3a and then gossip because of 4a.  On the other hand,
   if agent $i$ has not received a request, then agent $i$ must gossip virtually because of 4b.
Thus if $x_i(t) =  x_{i^*(t)}(t)$
   for some $i\in\scr{J}$, either a gossip or virtual gossip
   will have taken place between two neighboring agents with  an edge in
   $\bbb{G}$.
   %namely agents $i$ and $i^*(t)$.
   To complete the proof it
 is  thus enough  to consider the case when  $x_i(t)\neq x_{i^*(t)}(t)$ for all  $i\in\scr{J}$.
 %To prove the lemma is sufficient to show that a gossip must take place between
 % a pair $(i,i^*(t))$  with $i\in\scr{J}$.
We claim that under this condition at least one agent
with label $i\in\scr{J}$,
must place a request to gossip.
To prove that this is so, suppose the contrary.
Then there is  no agent with a label in $\scr{J}$ which is a requester so
$x_{i}(t)<x_{i^*(t)}(t)$ for all $i\in\scr{J}$. In particular
 $x_{i_1}(t)<x_{i_2}(t)<\ldots <
 x_{i_k}(t)<x_{i_k^*}$ where $\{i_1,i_2,\ldots,i_w\} = [i]_t$ and $k$ is the largest integer
  greater than $1$ for which  the labels $i_1,i_2,\ldots, i_{k}$ are all in $\scr{J}$.
Since the labels in $[i]_t$ are  all in $\scr{J}$, it must be that $k=w$  so $x_{i_1}(t)<x_{i_2}(t)<\ldots <
 x_{i_w}(t)<x_{i_w^*}(t)$. But $i_w^*(t)$ must equal some integer $i_j\in\{i_1,i_2,\ldots,i_{w-1}\}$
   so $x_{i_w}(t)<x_{i_j}(t)$.
   This
     is impossible because
 $j<w$. Therefore at least one agent  with a label in $\scr{J}$  must place a
 request to gossip.

To complete the proof it is enough to show that among the agents
 with labels in $\scr{J}$
 who receive requests to gossip at time $t$,
at least one agent - say agent $k$ - does not place a request to gossip.
 For if agent $k$ does not place a
request, then agent $k$ must gossip with that  agent
 with label closest to the front of $\mathbf{q}_k(t)$  who
  placed a request to gossip with agent $k$ at time $t$.

To
  prove that at least one agent receiving a gossip request at time $t$ does not place a request
   to gossip at time $t$,
   assume the contrary. Therefore suppose that every agent receiving a request to gossip a time $t$, also places
    a request to gossip at time $t$.
    Let $i$ be the label of any agent receiving a request to gossip at time $t$ and let
    $\{i_1,i_2,\ldots,i_w\} = [i]_t$.
Since agent $i_1 = i $  and $i$ receives a
 request to gossip, it also must place a request to gossip.  Hence agent $i_2$ must receive a request to gossip.
 Therefore agent $i_2$ must place a request to gossip at time $t$.
 %Therefore assuming agent $i_2$ must place a request to gossip.
By this reasoning one concludes that all of the agents  with labels $i_1,i_2,\ldots, i_w$
place requests to gossip at time $t$.  This implies that $x_{i_1}(t)>x_{i_2}(t)>\ldots >
 x_{i_w}(t)>x_{i_w^*(t)}(t)$.
 But $i_w^*(t)$ must equal some integer $i_j\in\{i_1,i_2,\ldots,i_{w-1}\}$. This means that
  $x_{i_w}(t)>x_{i_j}(t)$  with is impossible because
 $j<w$. Therefore at least one agent which has received a request to gossip has not placed a request to gossip.
\hfill
$\qed $

 It is worth noting that if $\bbb{G}$ has
 $s$ connected components, each with positive minimum degree, then there must be an edge $(a_i,b_i)$  within each
 component  for which  neighboring agents $a_i $ and $b_i$ either gossip or gossip virtually at time $t$.
  This can be proved
 using  an argument similar to the argument use to prove Lemma \ref{lib}.

\noindent{\bf Proof of Lemma \ref{pizza}:}
  We claim that $\bbb{A}$ satisfies
 the hypotheses of Lemma \ref{lib}.
Note first that
by assumption $\bbb{A}$ is a
connected graph with at least two vertices. Thus $\bbb{A}$
has at least one edge. Next observe that when
   $\bbb{G} = \bbb{A}$, we have
 $\scr{M}_i =\scr{N}_i$, $i\in\{1,2,\ldots,n\}$. Clearly $\bbb{A}$ automatically satisfies
 hypotheses of Lemma \ref{lib}. Hence Lemma \ref{pizza} is true.
\hfill
$\qed$

\noindent{\bf Proof of Proposition \ref{reeds}:}
See the proof of Proposition 2 in \cite{cdc11req}.
\hfill
$\qed $

Both analytical results and computer studies show that a slightly more complicated
queue update rule can lead to significantly faster convergence.

%Theorems  \ref{tm11} and \ref{ji}   thus imply that every sequence of gossip vectors generated
%by  Protocol III  converges to the desired limit point exponentially fast at a rate no worse that some finite
% number $\lambda <1$ which depends only $\bbb{A}$.  Calculation of this worst case
% bound is a subject for future research.

\subsubsection{Convergence Rate}

Theorems  \ref{tm11} and \ref{ji}  imply that every sequence of gossip vectors generated
by  Protocol III  converges to the desired limit point exponentially fast at a rate
no worse that some finite  number $\lambda <1$ which depends only on $\bbb{A}$.
%By Theorem 9 of \cite{pieee}, it is easy to see that
%$$\lambda=\max_{S\in\scr{C}}\sqrt[m]{\mu\{S\}}$$
%where $m$ is the number of edges of $\bbb{A}$, $\scr{C}$ is the compact set of
% complete gossip matrices determined by any gossiping sequence which has occurred over
% an interval of length $m$, and for any square matrix $T$, $\mu\{T\}$ denotes the second
%largest singular value of $T$. Thus one possible way to derive a worst case bound of
%$\lambda$ is to find a uniform bound of the second largest singular value
%of any complete gossip matrix in $\scr{C}$. In the sequel we will take another approach which
%appeals to a discrete-time Lyapunov function.
In the sequel we will derive a worst case bound of $\lambda$.

It is useful to think of a gossiping process in geometric terms.
Associate with agent $i$'s current gossip variable $x_i$, a corresponding point $x_i$ on the real line
which we will henceforth refer to as agent $i$'s current position. For agents $i$ and $j$ to gossip then
means simply that each moves to the midpoint between the two.

We would like to have a way to keep track of the entire group's progress in
reaching a consensus. Towards this end let us agree to call a nonnegative valued
function $V: \R^n\rightarrow\R$, an {\em indicator} if $V(t)=0$ just in case all
agents are at the same position at time $t$. In the sequel we will be concerned exclusively
with indicators comprised of sums of distances between pairs of points, and for
now we will assume that the specific pairs of points in question do {\em not}
change with iterations. Although we will talk exclusively about such functions
it is useful here to make this precise. So let $\scr{E}$ be a given subset of
$\{1,2,\ldots,n\}\times\{1,2,\ldots,n\}$; let us agree to say that a function
$V: \R^n\rightarrow [0,\infty)$ of the form
$$V(t)=\sum_{(i,j)\in\scr{E}}|x_i(t)-x_j(t)|$$
is a {\em multi-distance indicator function} if $V(t)=0$ implies that
all the $x_i$ have the same value at time $t$.
In the sequel we will often drop the modifier ``multi-distance''.
There is a natural way to associate with any such an indicator a simple \{undirected\}
graph. Specifically the {\em graph} of $V$, written $\bbb{G}_V$, is a that graph
on $n$ vertices, which has an edge $(i,j)$ just in case the distance between points
$i$ and $j$ is one of the terms in the sum comprising $V$. In other words,
$\scr{E}$ is the edge set of $V$.
Our first result characterizes the type of indicator functions under discussion.

\begin{lemma}
Let $\scr{E}\subset\{1,2,\ldots,n\}\times\{1,2,\ldots,n\}$ be fixed. A necessary and sufficient
condition for the function
\eq{V(t)=\sum_{(i,j)\in\scr{E}}|x_i(t)-x_j(t)|\label{sum}}
to be an indicator function is that $\bbb{G}_V$ contain a spanning tree of $\bbb{A}$.
\label{iff}\end{lemma}

\noindent
{\bf Proof of Lemma \ref{iff}:}
Suppose that $V$ is given by \rep{sum} and that $\bbb{G}_V$ contains a spanning tree $\bbb{T}\subset\bbb{A}$.
If $V (t) = 0$ then $|x_i(t)-x_j(t)| = 0$, $(i,j)\in\scr{E}$. Since $\bbb{G}_V$ contains a spanning tree of
$\bbb{A}$ and $\bbb{A}$ is connected, this can only occur if all agents are in the same position.
Consequently $x = 0$ and $V$ is an indicator.

For the converse, suppose that $V$ is an indicator and that $\bbb{G}_V$ does not contain a spanning tree
of $\bbb{A}$. This means that $\bbb{G}_V$ must contain a connected component whose vertex set $\scr{U}$ is disjoint
from the vertex set of the union of the vertex sets of the remaining connected components of $\bbb{G}_V$.
Moreover, there must be an edge $(i,j)$ of $\bbb{A}$ such that $i\in\scr{U}$ and $j\notin\scr{U}$. Now pick two distinct
points $y_1$ and $y_2$ and position all agents with indices in $\scr{U}$ at $y_1$ and all agents with indices not in
$\scr{U}$ at $y_2$. Then $V (t) = 0$ because $|x_i-x_j|$ is not in the sum defining $V$. But all agents are not in
the same position which is a contradiction of the hypothesis that $V$ is an indicator. Therefore $\bbb{G}_V$
must contain a spanning tree of $\bbb{A}$.
\hfill
$\qed$

%Let $x$ be the agents' current position vector and $y$ the agents' position vector
%after agents $i$ and $j$ have moved.
Suppose agents $i$ and $j$ gossip at time $t$.
Let us say that an indicator $V$ is
{\em instantaneous} if there is a positive number $\lambda$ such that
\eq{V(t+1)-V(t)\leq -\lambda|x_i(t)-x_j(t)|\label{V}}
Thus if $V$ is instantaneous, there is a definite decrease in its value whenever any
allowable pair of agents not initially in the same position, gossip.
In the sequel it will be shown that if $V$ is instantaneous, then $\lambda$ can
be taken as $1$ for all allowable gossips.

\begin{lemma}
If $V$ is instantaneous, then $\bbb{A}\subset\bbb{G}_V$.
\label{contain}\end{lemma}

\noindent
{\bf Proof of Lemma \ref{contain}:}
Let $V$ be instantaneous. Then it is easy to see that  $\bbb{G}_V$ must be connected.
Suppose that $(i,j)$ is an edge in $\bbb{A}$ which is not in $\bbb{G}_V$.
Position agents $i$ and $j$ so that $x_j(t) > x_i(t)$.
Let $\scr{K}$ be the set of vertices $k$ such that both $(i,k)$ and $(j,k)$ are edges in $\bbb{G}_V$.
For $k\in\scr{K}$, position agent $k$ at a point between $x_i$ and $x_j$.
For $k\notin\scr{K}$, position agent $k$ at a point $x_k > x_j$ if $(j,k)$ is an
edge in $\bbb{G}_V$ or at a point $x_k < x_i$ if $(i,k)$ is an edge in $\bbb{G}_V$.
Position all remaining agents at any fixed points. Because $\bbb{G}_V$ is connected,
this positions all $n$ agents. Note that after the gossip, for
all $k\in\scr{K}$, the sum of the distances from agent $k$ to agents $i$ and $j$ is the same as before the gossip.
This is also true of any agent which is not a neighbor of agent $i$ or $j$. Meanwhile, for $k\notin\scr{K}$, the
distances from agent $k$ to agent $j$ increases if $(j,k)$ is an edge in $\bbb{G}_V$ as does the distance from agent
$k$ to agent $i$ if $(i,k)$ is an edge in $\bbb{G}_V$. Since $(i,j)$ is not an edge in $\bbb{G}_V$,
this means that the distance
between each two agents which are neighbors in $\bbb{G}_V$ does not decrease.
Therefore $V(t+1)-V(t)\geq 0$.
Meanwhile the distance before gossip between agents $i$ and $j$ is positive. Therefore there is no
positive number $\lambda$ for which \rep{V} holds. Thus $V$ is not instantaneous which is a contradiction. Thus
the lemma is true.
\hfill
$\qed$

\begin{proposition}
If $V$ is an instantaneous indicator and $i$ and $j$ are an allowable pair of agents
who gossip at time $t$, then
\eq{V(t+1)-V(t)\leq -|x_i(t)-x_j(t)|\label{V1}}
\label{single}\end{proposition}
The proof of this proposition depends on the following result.

\begin{lemma}
Suppose that agents $i$ and $j$ gossip. Let $k$ be different than $i$ and $j$. Then the distance
between $k$ and $i$ plus the distance between $k$ and $j$ after the gossip is no greater than the distance
between $k$ and $i$ plus the distance between $k$ and $j$ before the gossip.
\label{yu1}\end{lemma}

\noindent
{\bf Proof of Lemma \ref{yu1}:}
First consider the case when before gossip, agent $k$ is not in between agents $i$
and $j$. Without loss of generality assume that $x_i \leq x_j$. If $x_k \leq x_i$, then before gossip, the distance
sum is $x_i-x_k + x_j-x_k = x_i + x_j - 2x_k$ whereas after the gossip the sum is
$2\frac{1}{2} (x_i - x_j) - 2x_k = x_i+x_j -2x_k$. If $x_j \leq x_k$,
then before gossip, the distance sum is $x_k -x_i+x_k -x_j = 2x_k -(x_i+x_j)$
whereas after the gossip the sum is $2x_k - 2\frac{1}{2} (x_i - x_j) = 2x_k - (x_i + x_j)$. Therefore in either case
gossiping does not affect the distance sum.

Now consider the case when agent $k$ is initially in between agents $i$ and $j$. Without loss of
generality assume that $x_i \leq x_k \leq x_j$. Then before gossiping, the distance sum is
$x_j -x_k+x_k-x_i = x_j - x_i$ whereas after the gossip, the sum is
$2\frac{1}{2} (x_i + x_j) - 2x_k = x_i + x_j - 2x_k$. But $x_i \leq x_k$ so
$x_i + x_j - 2x_k \leq x_j - x_i$.
\hfill
$\qed$

\noindent
{\bf Proof of Proposition \ref{single}:}
Suppose that agents $i$ and $j$ gossip. In view of Lemma \ref{yu1}, the sum total
of all distances appearing in $V$, with the exception of the distance between agents $i$ and $j$, does not
increase after the gossip. Meanwhile the distance between agents $i$ and $j$ decreases by $|x_i - x_j|$.
Since $\bbb{A}\subset \bbb{G}_V$ \{by Lemma \ref{contain}\},
the distance between $i$ and $j$ must also appear in the definition of $V$. Therefore \rep{V1} must hold.
\hfill
$\qed$

\begin{proposition}
A necessary and sufficient condition for $V$ to be an instantaneous indicator is that
$\bbb{A}\subset\bbb{G}_V$ and for each edge $(i,k)$ of $\bbb{G}_V$ for which $(i,j)$
is an edge of $\bbb{A}$, $(j,k)$ is an edge of $\bbb{G}_V$.
\label{morse}\end{proposition}

\noindent
{\bf Proof of Proposition \ref{morse}:}
Suppose that $V$ is an indicator with the properties that $\bbb{A}\subset\bbb{G}_V$ and
each edge $(i,k)$ of $\bbb{G}_V$ for which $(i,j)$ is an edge of $\bbb{A}$, $(j,k)$ is an
edge of $\bbb{G}_V$. Suppose that agents $i$ and $j$ gossip in which case $(i,j)$ is an
edge of $\bbb{A}$ and thus $\bbb{G}_V$. Let $(m,k)$ be any edge in $\bbb{G}_V$.
If $\{i,j\}$ and $\{m,k\}$ are disjoint sets, the distance between agents $m$ and $k$
does not change with the gossip. If $\{i,j\}$ and $\{m,k\}$ are not disjoint sets,
then without loss of generality we can take $m=i$. Thus by hypothesis both
$(i,k)$ and $(j,k)$ are edges in $\bbb{G}_V$. But by Lemma \ref{yu1} the sum of the distance
between agent $k$ and agent $i$ and the distance between agent $k$ and agent $j$ does not increase after the
gossip. Since this is true for all edges in $\bbb{G}_V$ with the exception of $(i,j)$, it must be true that
\rep{V} holds with $\lambda = 1$. Therefore $V$ is instantaneous.
The simple proof of the necessity part of this proposition is omitted.
%Now suppose that $V$ is an instantaneous indicator.
%Then $\bbb{A}\subset\bbb{G}_V$ because of Lemma \ref{contain}.
\hfill
$\qed$

\begin{theorem}
$V$ is instantaneous if and only if $\bbb{G}_V$ is complete.
\label{anu}\end{theorem}

\noindent
{\bf Proof of Theorem \ref{anu}:}
Suppose that $(i,j)$ is not an edge in $\bbb{G}_V$. Since $\bbb{A}$ is connected, there
must be a path from $i$ to $j$ in $\bbb{A}$ and thus $\bbb{G}_V$.
Suppose there are other $k>0$ vertices in the path. Then the path consists of $k+1$ edges
which are denoted by $(i, v_1), (v_1,v_2),\ldots, (v_{k-1},v_k),(v_k, j)$.
%The $k+1$ edges are all contained in $\bbb{A}$ and $\bbb{G}_V$.
Since $(i,v_1)$ is in $\bbb{G}_V$ and $(v_1,v_2)$ is in $\bbb{A}$, then by Proposition \ref{morse},
$(i,v_2)$ is in $\bbb{G}_V$. Similarly, since $(i,v_2)$ is in $\bbb{G}_V$ and
$(v_2,v_3)$ is in $\bbb{A}$, then $(i,v_3)$ is also in $\bbb{G}_V$.
By repeating this argument, one reaches the conclusion that $(i,j)$ is an edge of $\bbb{G}_V$,
which is a contradiction. Thus  $\bbb{G}_V$ must be a complete graph.
\hfill
$\qed$

By Theorem \ref{anu}, it is clear that the desired instantaneous indicator must
be in the form of
$$V(t)=\sum_{(i,j)\in\scr{A}}|x_i(t)-x_j(t)|$$
where $\scr{A}=\{1,2,\ldots,n\}\times\{1,2,\ldots,n\}$.

\begin{lemma}
Suppose all $n$ agents follow Protocol III. Let $m$ be the number of edges of $\bbb{A}$.
Then for any time $t$,
$$V(t+m)\leq \displaystyle\left(1-\frac{4}{n^2}\right)V(t)$$
\end{lemma}
A proof of this lemma can be found in \cite{low&murray}
(see Lemma~2 in \cite{low&murray}).

We are led to the following result.

\begin{proposition}
Suppose all $n$ agents follow Protocol III.
Then every sequence of gossip vectors $x(t)$, $t>0$ generated  converges
to the desired limit point exponentially fast at a rate no worse than
$$\displaystyle\left(1-\frac{4}{n^2}\right)^{\frac{1}{m}}$$
where $m$ is the number of edges of $\bbb{A}$.
\end{proposition}

\section{Concluding Remarks}\label{end}

Three request-based gossiping protocols with different types of
queue updates are studied, which provides an in-depth understanding of how
local ordering and queue updates avoid deadlocks.
It is shown that the rule of queue updates has significant effects on
convergence and convergence time.
One of the problems with the idea of gossiping, which apparently is not widely appreciated,
 is that it is difficult to devise  provably correct
gossiping protocols  which are guaranteed to avoid deadlocks without making restrictive assumptions.
 The research in this paper and in
 \cite{low&murray,alex} contributes to our understanding of this issue and how to deal with it.
For the protocols presented in this paper, it is assumed that the communication between
agents is delay-free. Analysis of the effect of transmission delays is a subject for future research.

\section{Acknowledgement}

The authors wish to thank Ming Cao (University of Groningen)
for useful discussions which have contributed to this work.

\bibliographystyle{unsrt}
\bibliography{pieee}

\end{document}